\documentclass[11pt]{amsart}%
\usepackage{graphicx,amscd,color,amsmath,amsfonts,amssymb,geometry}
\usepackage{amsmath}
\usepackage{amsfonts}
\usepackage{amssymb}
\usepackage{graphicx}%
\providecommand{\U}[1]{\protect\rule{.1in}{.1in}}
\newtheorem{theorem}{Theorem}[section]

\theoremstyle{definition}

\newtheorem{remark}[theorem]{Remark}

\thanks{D. Pellegrino is supported by CNPq Grant 307327/2017-5 and Grant 2019/0014 Paraiba State Research Foundation (FAPESQ)}
\thanks{J. Silva is supported by CAPES}
\begin{document}
\title[On a continuous Gale--Berlekamp switching game]{On a continuous Gale--Berlekamp switching game}
\author[Daniel\ Pellegrino, Janiely Silva and Eduardo V. Teixeira]{Daniel\ Pellegrino, Janiely Silva and Eduardo V. Teixeira}
\address{D. Pellegrino and J. Silva\\
Departamento de Matem\'atica\\
Universidade Federal da Para\'{\i}ba \\
Brazil\\
E. V. Teixeira\\
Department of Mathematics\\
University of Central Florida\\
Estados Unidos.}
\email{pellegrino@pq.cnpq.br}
\email{janiely.silva@estudantes.ufpb.br}
\email{Eduardo.Teixeira@ucf.edu}

\begin{abstract}
We propose a continuous version of the classical Gale--Berlekamp switching
game. We also study a weighted version of this new continuous game. The main results of this paper concern growth estimates for the corresponding optimization problems. The methods developed in this article  are deterministic in nature and in some special cases the estimates obtained are optimal. 

\end{abstract}
\maketitle
\tableofcontents

\section{Introduction}

Designed independently by Elwyn Berlekamp and David Gale in the $1960$'s, the
Gale--Berlekamp switching game -- also known as the unbalancing lights problem
-- represents a trademark in the field of combinatorics and its applications,
with special interests in the theory of Computer Sciences. This single-player
game consists of an $n\times n$ square matrix of light bulbs set-up at an
initial light configuration. The goal is to turn off as many lights as
possible using $n$ row and $n$ column switches, which invert the state of each
bulb in the corresponding row or column.

For an initial pattern of lights $\Theta$, let $i(\Theta)$ denote the smallest
final number of on-lights achievable by row and column switches starting from
$\Theta.$ The smallest possible number of remaining on-lights $R_{n}$,
starting from a worse initial pattern, is then
\[
R_{n}=\max\{i(\Theta):\Theta\text{ is an }n\times n\text{ light pattern}\}.
\]
Sometimes this optimization problem is posed as to find the maximum of the
difference between lights on and off, often denoted by $G_{n}$. Obviously both
problems are equivalent as $R_{n}=\frac{1}{2}\left(  n^{2}-G_{n}\right)  $.

The original problem introduced by Berlekamp asks for the exact value of
$R_{10}$ and it was proved in \cite{jor} that $R_{10}=35$ (and thus
$G_{10}=30$)$.$ Several related questions pertaining to the original problem
have been investigated in depth, see e.g. \cite{br, jor, fish, sch}, in
particular the hardness of solving the Gale-Berlekamp game was studied in
\cite{roth}.

In this paper we propose a continuous version of the Elwyn Berlekamp switching
game. We are interested in a continuous version of the game for which vectors replace
light bulbs and direction knobs substitute the
discrete switches, used to invert the state of the bulbs in the
original problem. We also allow for non-square game-boards.

To explain the new proposed game, we initially notice that by
associating $+1$ to the on-lights and $-1$ to the off-lights from the array of
lights $\left(  a_{ij}\right)  _{i,j=1}^{n}$ the original game can be treated
mathematically as
\[
G_{n}=\min\left\{  \max_{x_{i},y_{j}\in\{-1,1\}}\left\vert \sum\limits_{i,j=1}%
^{n}a_{ij}x_{i}y_{j}\right\vert :a_{ij}=-1\text{ or } +1\right\}  ,
\]
where $x_{i}$ and $y_{j}$ denote the switches of row $i$ and of column $j$, respectively.

The new optimization problem herein proposed involves a matrix $\left(
a_{ij}\right)  $ with $n_{1}$ rows and $n_{2}$ columns whose elements are unit
vectors of the plane, $\mathbb{R}^{2}.$ The initial direction pattern of each
$n_{1}n_{2}$ vectors is set up at the beginning of the game. For each row $i$
and each column $j$ there are knobs $x_{i}$ and $y_{j}$, respectively.
Rotating knob $x_{i}$ by an angle $\theta_{i},$ it affects all vectors
$a_{ij}$ of the row $i$. Analogously, when knob $y_{j}$ is rotated by an angle
$\theta_{j},$ the same happens with all the vectors $a_{ij}$ of the column $j$
(see Figure 1.).

\begin{figure}[tbh]
\vspace*{0,2cm} \centering
\includegraphics[width=0.5\textwidth]{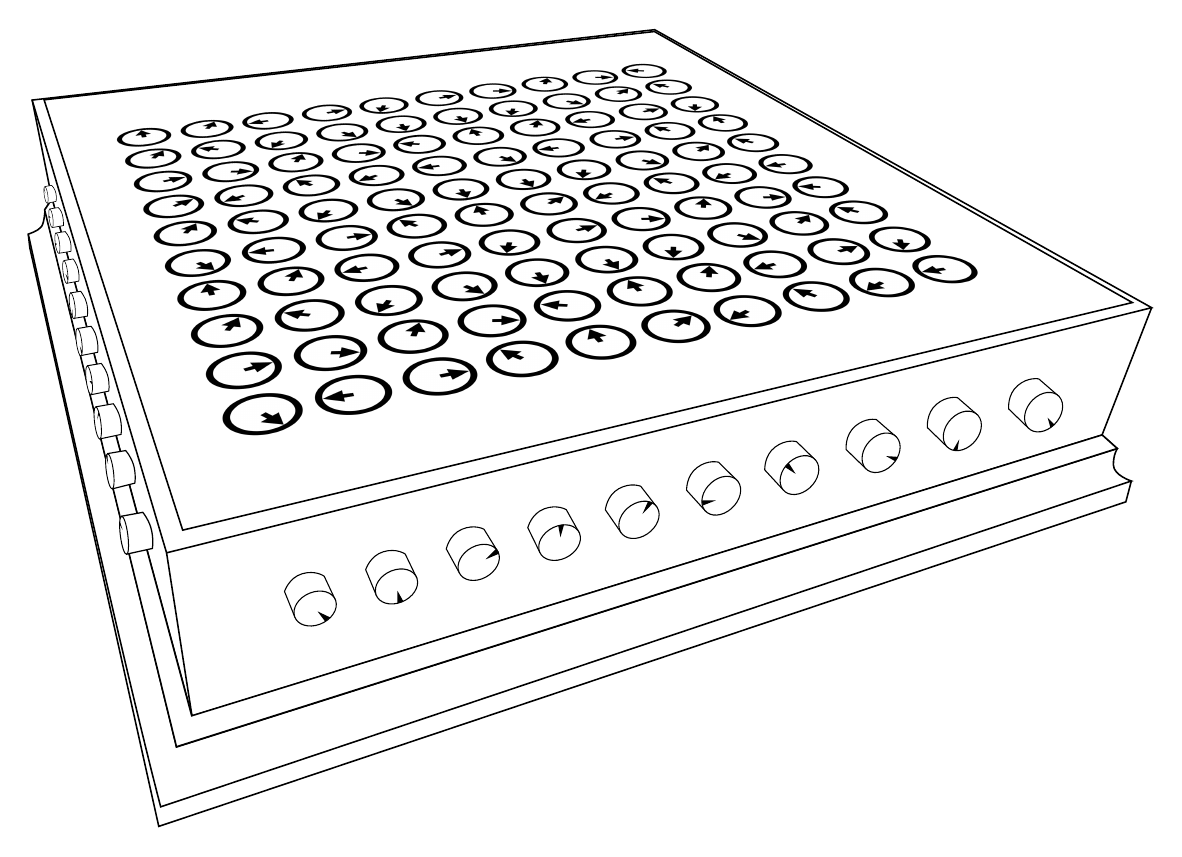} \medskip
{\small \centerline{Figure 1.} }\end{figure}\bigskip

The game consists of maximizing the Euclidean norm of the sum of all vectors
in the final stage. More precisely, for an initial pattern $\Theta$ of unit
vectors, let $s(\Theta)$ be the supremum of the (Euclidean) norms of the sums
of all $n_{1}n_{2}$ vectors achievable by row and column adjusts. The extremal
problem, that is if one starts with the worse possible initial pattern, is to
determine
\[
G_{n_{1}n_{2}}^{(1)}:=\min\{s(\Theta):\Theta\text{ an }n_{1}\times n_{2}\text{
pattern}\}.
\]

Regarding this problem, our main result provides an asymptotic growth for
$G_{n_{1}n_{2}}^{(1)}$, as indicated in the following theorem:

\begin{theorem}
\label{6b}For all integers $n_{1}, n_{2}$, we have
\begin{equation}
0.886\leq\frac{G_{n_{1}n_{2}}^{(1)}}{\sqrt{n_{1}n_{2}}\max\{\sqrt{n_{1}}%
,\sqrt{n_{2}}\}}\leq1. \label{x1}%
\end{equation}

\end{theorem}

In this paper we are also interested in a {weighted version} of
the above problem. More precisely, we want to understand growth estimates when {the knobs in the $n_{2}$ columns may be used not only to rotate column $j$ ,  but also multiply the vectors  $a_{ij}$ by a choice of real numbers $b_{j}$} verifying
\begin{equation*}
{\textstyle\sum\limits_{j=1}^{n_{2}}}
b_{j}^{2}=1.\label{pp}%
\end{equation*}

Again the game consists of maximizing the Euclidean norm of the sum of all
vectors. In other words, for an initial pattern $\Theta$ of unit vectors, let
$s(\Theta)$ be the supremum of the Euclidean norms of the sums of all
$n_{1}n_{2}$ vectors achievable by row and column adjusts. The extremal
problem, starting at the worse possible initial pattern, is to determine
\[
G_{n_{1}n_{2}}^{(2)}:=\min\{s(\Theta):\Theta\text{ an }n_{1}\times n_{2}\text{
pattern}\}.
\]

The {weighted} feature of the second game yields
 {a different freedom to the player}. Surprisingly, in this case, we shall obtain a definitive, constructive solution to the problem. Indeed, the second main result we prove in this article is:

\begin{theorem}
\label{7b}For all positive integers $n_{1},n_{2}\geq2$, with $n_{2}\geq n_{1}
$, we have
\begin{equation}
\frac{G_{n_{1}n_{2}}^{(2)}}{\sqrt{n_{1}n_{2}}}=1. \label{x2}%
\end{equation}
\end{theorem}

We conclude this introduction by commenting on the ideas and techniques used
to prove both Theorems \ref{6b} and \ref{7b}, which are of particular
interest. We observe  due to the combinatorial complexity of this kind of problems, growth estimates as in Theorem \ref{6b} and Theorem \ref{7b} are often obtained by non-deterministic techniques, see for instance \cite{alon,ap5,ben2}.
A main novelty proposed in this article regards a deterministic approach to
estimate both $G_{n_{1}n_{2}}^{(1)}$ and $G_{n_{1}n_{2}}^{(2)}$, which yields
improved, more precise estimates than non-deterministic methods. We believe methods herein developed are likely to be applicable in an array of other problems and to exemplify
the depth of these new ideas, we also prove analogues of (\ref{x1}) and
(\ref{x2}) in higher-dimensional spaces.

\section{Proof of Theorem \ref{6b}}

Initially, it is more convenient to conceive the vectors in the game as
complex numbers $a_{ij}$ with norm $1$, which represent the elements of the
array $\left(  a_{ij}\right)  _{i,j=1}^{n_{1},n_{2}}$. In this case, when the player rotates a knob, the action is modeled by multiplication by unimodular
complex numbers. In the  {weighted} game, to be treated in the
next section, the action of {rotating and multiplying} by real numbers $b_{j},$ with
${\textstyle\sum\limits_{j=1}^{n_{2}}} b_{j}^{2}=1$, will then be modeled as
multiplication by {certain} complex numbers with norm less than or equal to $1$.

We start off the proof of Theorem \ref{6b} by reminding that a consequence of
the Krein--Milman Theorem assures that for all $A \colon\ell_{\infty}^{n_{1}%
}\times\ell_{\infty}^{n_{2}}\rightarrow\mathbb{C}$, there holds%
\[
\left\Vert A\right\Vert =\sup_{\left\vert x_{j_{1}}^{(1)}\right\vert
=\left\vert x_{j_{2}}^{(2)}\right\vert =1}\left\vert \sum\limits_{j_{1}%
,j_{2}=1}^{n_{1},n_{2}}a_{j_{1}j_{2}}x_{j_{1}}^{(1)} x_{j_{2}}^{(2)}%
\right\vert
\]
and thus,
\[
G^{(1)}_{n_{1}n_{2}}=\inf\left\{  \left\Vert A\right\Vert :\left\vert
a_{j_{1}j_{2}}\right\vert =1\right\}  .
\]

Our task is then to estimate $\inf\left\{  \left\Vert A\right\Vert :\left\vert
a_{j_{1}j_{2}}\right\vert =1\right\}  $, where the infimum runs over all
bilinear forms \linebreak$A:\ell_{\infty}^{n_{1}}\times\ell_{\infty}^{n_{2}%
}\rightarrow\mathbb{C}$ with unimodular coefficients.

Once the problem has been described as above, the upper bound in Theorem
\ref{6b} can be obtained by means of an argument from the seminal paper of
Bohnenblust and Hille \cite[Theorem II, page 608]{bhv}. We shall explain the
necessary adaptations when we deliver the proof of Theorem \ref{8b}, on
Section \ref{Section Anisotropic}.

As for the lower estimate, we shall make use of Khinchin inequality, which we
revise for the sake of completeness.

\subsection{Khinchin inequality}

To motivate, let's state the following question: suppose that we have $n$ real
numbers $a_{1},\ldots,a_{n}$ and a fair coin. When we flip the coin, if it
comes up heads, you choose $\beta_{1}=a_{1}$, and if it comes up tails, you
choose $\beta_{1}=-a_{1}.$ When we play for the second time, if it comes up
heads, you choose $\beta_{2}=\beta_{1}+a_{2}$ and, if it comes up tails, you
choose $\beta_{2}=\beta_{1}-a_{2}$. Repeating the process, after having
flipped the coin $k$ times we have
\[
\beta_{k+1}:=\beta_{k}+a_{k+1},
\]
if it comes up heads and
\[
\beta_{k+1}:=\beta_{k}-a_{k+1},
\]
if it comes up tails. After $n$ steps, what should be the expected value of
\[
|\beta_{n}|=\left\vert \sum_{k=1}^{n}\pm a_{k}\right\vert ?
\]
Khinchin's inequality, see for instance
{\cite[page 10]{diestel}}, shows that the \textquotedblleft
average\textquotedblright\
\[
\left(  \frac{1}{2^{n}}\left\vert \sum\limits_{\varepsilon_{1},\ldots
,\varepsilon_{n}=1,-1}\varepsilon_{j}a_{j}\right\vert ^{p}\right)  ^{\frac
{1}{p}}
\]
behaves as the $\ell_{2}$-norm of $\left(  a_{n}\right)  .$ More precisely, it
asserts that for any $p>0$ there are constants $A_{p},B_{p}>0$ such that
\begin{equation}
A_{p}\left(  \sum\limits_{j=1}^{n}\left\vert a_{j}\right\vert ^{2}\right)
^{\frac{1}{2}}\leq\left(  \frac{1}{2^{n}}\left\vert \sum\limits_{\varepsilon
_{1},\ldots,\varepsilon_{n}=1,-1}\varepsilon_{j}a_{j}\right\vert ^{p}\right)
^{\frac{1}{p}}\leq B_{p}\left(  \sum\limits_{j=1}^{n}\left\vert a_{j}%
\right\vert ^{2}\right)  ^{\frac{1}{2}} \label{65}%
\end{equation}
for all sequence of scalars $\left(  a_{i}\right)  _{i=1}^{n}$ and all
positive integers $n.$

\bigskip

Back to the proof of Theorem \ref{6b}, for the purpose of establishing a lower
estimate for the growth of $G^{(1)}_{n_{1} n_{2}}$, we are interested in the
case $p=1.$ The natural counterpart for the average $\frac{1}{2^{n}}\left\vert
\sum\limits_{\varepsilon_{1},\ldots,\varepsilon_{n}=1,-1}\varepsilon_{j}%
a_{j}\right\vert $ in the complex framework is
\begin{equation}
\left(  \frac{1}{2\pi}\right)  ^{n}\int_{0}^{2\pi}\ldots\int_{0}^{2\pi
}\left\vert \sum_{j=1}^{n}a_{j}e^{it_{j}}\right\vert dt_{1}\cdots dt_{n}.
\label{222}%
\end{equation}
It is well known that in this new context we also have a Khinchin-type
inequality, called Khinchin inequality for Steinhaus variables, which asserts
that there exist constants $\widetilde{A_{1}}$ and $\widetilde{B_{1}}$ such
that
\begin{equation}
\widetilde{A_{1}}\left(  \sum_{j=1}^{n}\left\vert a_{j}\right\vert
^{2}\right)  ^{\frac{1}{2}}\leq\left(  \frac{1}{2\pi}\right)  ^{n}\int
_{0}^{2\pi}\ldots\int_{0}^{2\pi}\left\vert \sum_{j=1}^{n}a_{j}e^{it_{j}%
}\right\vert dt_{1}\cdots dt_{n}\leq\widetilde{B_{1}}\left(  \sum_{j=1}%
^{n}\left\vert a_{j}\right\vert ^{2}\right)  ^{\frac{1}{2}} \label{est}%
\end{equation}
for every positive integer $n$ and scalars $a_{1},\ldots,a_{n}$. In \cite{kon}
it is proven that $\widetilde{A_{1}}=\frac{\sqrt{\pi}}{2}$ and, for our
purposes, from now on we shall only be interested in the left hand side of
(\ref{est}). For a bilinear form $A:\ell_{\infty}^{n_{1}}\times\ell_{\infty}^{n_{2}%
}\rightarrow\mathbb{C}$ given by%
\[
A\left(  e_{j_{1}},e_{j_{2}}\right)  =a_{j_{1}j_{2}}%
\]
with%
\[
\left\vert a_{j_{1}j_{2}}\right\vert =1,
\]
we have
\begin{align*}
&  \left(  \sum_{j_{1}=1}^{n_{1}}\left\vert A\left(  e_{j_{1}},e_{j_{2}%
}\right)  \right\vert ^{2}\right)  ^{1/2}\\
&  \leq\left(  \frac{2}{\sqrt{\pi}}\right)  \left(  \frac{1}{2\pi}\right)
^{n_{1}}\int_{0}^{2\pi}\ldots\int_{0}^{2\pi}\left\vert \sum_{j_{1}=1}^{n_{1}%
}A\left(  e_{j_{1}},e_{j_{2}}\right) e^{it_{j_{1}}^{\left(  1\right)}}
\right\vert dt_{j_{1}}^{\left(  1\right)  }\\
&  =\left(  \frac{2}{\sqrt{\pi}}\right)  \left(  \frac{1}{2\pi}\right)
^{n_{1}}\int_{0}^{2\pi}\ldots\int_{0}^{2\pi}\left\vert A\left(  \sum_{j_{1}%
=1}^{n_{1}}e^{it_{j_{1}}^{\left(  1\right)  }}e_{j_{1}},e_{j_{2}}\right)
\right\vert dt_{j_{1}}^{\left(  1\right)  }.
\end{align*}
Since
\begin{align*}
&  \int_{0}^{2\pi}\ldots\int_{0}^{2\pi}\sum\limits_{j_{2}=1}^{n_{2}}\left\vert
A\left(  \sum_{j_{1}=1}^{n_{1}}e^{it_{j_{1}}^{\left(  1\right)  }}e_{j_{1}%
},e_{j_{2}}\right)  \right\vert dt_{j_{1}}^{\left(  1\right)  }\\
&  \leq\left(  2\pi\right)  ^{n_{1}}\max_{t_{j_{1}}^{\left(  1\right)  }%
\in\lbrack0,2\pi]}\sum\limits_{j_{2}=1}^{n_{2}}\left\vert A\left(  \sum
_{j_{1}=1}^{n_{1}}e^{it_{j_{1}}^{\left(  1\right)  }}e_{j_{1}},e_{j_{2}%
}\right)  \right\vert ,
\end{align*}
denoting the topological dual of ${\ell_{\infty}^{n}}$ by $\left(
{\ell_{\infty}^{n}}\right)  ^{\ast}$ and its closed unit ball by $B_{\left(
{\ell_{\infty}^{n}}\right)  ^{\ast}},$ we have
\begin{align*}
&  \sum\limits_{j_{2}=1}^{n_{2}}\left(  \sum\limits_{j_{1}=1}^{n_{1}%
}\left\vert A\left(  e_{j_{1}},e_{j_{2}}\right)  \right\vert ^{2}\right)
^{1/2}\\
&  \leq\left(  \frac{2}{\sqrt{\pi}}\right)  \left(  \frac{1}{2\pi}\right)
^{n_{1}}\int_{0}^{2\pi}\ldots\int_{0}^{2\pi}\sum\limits_{j_{2}=1}^{n_{2}%
}\left\vert A\left(  \sum_{j_{1}=1}^{n_{1}}e^{it_{j_{1}}^{\left(  1\right)  }%
}e_{j_{1}},e_{j_{2}}\right)  \right\vert dt_{j_{1}}^{\left(  1\right)  }\\
&  \leq\left(  \frac{2}{\sqrt{\pi}}\right)  \max_{t_{j_{1}}^{\left(  1\right)
}\in\lbrack0,2\pi]}\sum\limits_{j_{2}=1}^{n_{2}}\left\vert A\left(
\sum_{j_{1}=1}^{n_{1}}e^{it_{j_{1}}^{\left(  1\right)  }}e_{j_{1}},e_{j_{2}%
}\right)  \right\vert \\
&  \leq\left(  \frac{2}{\sqrt{\pi}}\right)  \left\Vert A\right\Vert
\sup_{\varphi\in B_{\left(  {\ell_{\infty}^{n}}\right)  ^{\ast}}}%
\sum\limits_{j_{2}=1}^{n_{2}}\left\vert \varphi\left(  e_{j_{2}}\right)
\right\vert \\
&  =\left(  \frac{2}{\sqrt{\pi}}\right)  \left\Vert A\right\Vert ,
\end{align*}
where in the last equality we have used the isometric isometry
\[%
\begin{array}
[c]{cccc}
& \ell_{1}^{n} & \longrightarrow & \left(  {\ell_{\infty}^{n}}\right)  ^{\ast
}\\
& \left(  a_{j}\right)  _{j=1}^{n} & \longmapsto & \varphi,
\end{array}
\]
with $\varphi:{\ell_{\infty}^{n}}\rightarrow\mathbb{C}$ defined by
\[
\varphi\left(  \left(  x_{j}\right)  _{j=1}^{n}\right)  =%
{\textstyle\sum\limits_{j=1}^{n}}
a_{j}x_{j}.
\]
Finally, since $\left\vert A\left(  e_{j_{1}},e_{j_{2}}\right)  \right\vert
=1$, we conclude that
\[
\left\Vert A\right\Vert \geq\left(  \frac{\sqrt{\pi}}{2}\right)  n_{2}%
n_{1}^{\frac{1}{2}}.
\]
Thus%
\[
\left(  \frac{\sqrt{\pi}}{2}\right)  \leq\frac{G_{n_{1}n_{2}}^{(1)}}%
{\sqrt{n_{1}n_{2}}\max\{\sqrt{n_{1}},\sqrt{n_{2}}\}}.
\]

\section{Proof of Theorem \ref{7b}}

For the upper estimate, consider an $n_{2}\times n_{2}$ matrix $(a_{rs})$ with
$a_{rs}=e^{2\pi i\frac{rs}{n_{2}}}$. Note that $|a_{rs}|=1$ and
\begin{equation}
\sum_{s=1}^{n_{2}}a_{rs}\overline{a_{ts}}=n_{2}\delta_{rt}. \label{delt}%
\end{equation}
Define
\[
A:\ell_{\infty}^{n_{1}}\times\ell_{2}^{n_{2}}\longrightarrow\mathbb{C}%
\]
by
\[
A\left( x^{(1)},x^{(2)}\right) =\sum_{i{_{1}},i_{2}=1}^{n_{1},n_{2}}%
a_{i_{1},i_{2}}x_{i_{1}}^{(1)}x_{i_{2}}^{(2)}.
\]
Let $x^{(1)}\in B_{\ell_{\infty}^{n_{1}}}$ and $x^{(2)}\in B_{\ell_{2}^{n_{2}%
}}$, then $y^{(1)}\in B_{\ell_{\infty}^{n_{2}}}$ and $y^{(2)}\in B_{\ell
_{2}^{n_{2}}}$, where
\[
y^{(1)}=\left(  x_{1}^{(1)},\ldots,x_{n_{1}}^{(1)},0,\ldots,0\right)
~~\mbox{and}~~y^{(2)}=x^{(2)}=\left(  x_{1}^{(2)},\ldots,x_{n_{2}}%
^{(2)}\right)  .
\]
Thus
\begin{align*}
\left\vert A\left(  x^{(1)},x^{(2)}\right)  \right\vert  &  =\left\vert
\sum_{i_{1},i_{2}}^{n_{2}}a_{i_{1}i_{2}}y_{i_{1}}^{(1)}y_{i_{2}}%
^{(2)}\right\vert \\
&  \leq\sum_{i_{2}=1}^{n_{2}}\left\vert \sum_{i_{1}=1}^{n_{2}}a_{i_{1}i_{2}%
}y_{i_{1}}^{(1)}\right\vert \left\vert y_{i_{2}}^{(2)}\right\vert \\
&  \leq\left(  \sum_{i_{2}=1}^{n_{2}}|y_{i_{2}}^{(2)}|^{2}\right)  ^{\frac
{1}{2}}\cdot\left(  \sum_{i_{2}=1}^{n_{2}}\left\vert \sum_{i_{1}=1}^{n_{2}%
}a_{i_{1}i_{2}}y_{i_{1}}^{(1)}\right\vert ^{2}\right)  ^{\frac{1}{2}}\\
&  \leq\left(  \sum_{i_{2}=1}^{n_{2}}\left\vert \sum_{i_{1}=1}^{n_{2}}%
a_{i_{1}i_{2}}y_{i_{1}}^{(1)}\right\vert ^{2}\right)  ^{\frac{1}{2}}.
\end{align*}
Since%

\[
\left(  \sum_{i_{2}=1}^{n_{2}}\left\vert \sum_{i_{1}=1}^{n_{2}}a_{i_{1}i_{2}%
}y_{i_{1}}^{(1)}\right\vert ^{2}\right)  ^{\frac{1}{2}}=\left(  \sum_{i_{2}%
=1}^{n_{2}}\sum_{\substack{i_{1}=1\\j_{1}=1}}^{n_{2}}y_{i_{1}}^{(1)}%
\overline{y_{j_{1}}^{(1)}}a_{i_{1}i_{2}}\overline{a_{j_{1}i_{2}}}\right)
^{\frac{1}{2}}=\left(  \sum_{\substack{i_{1}=1\\j_{1}=1}}^{n_{2}}y_{i_{1}}%
^{1}\overline{y_{j_{1}}^{(1)}}\sum_{i_{2}=1}^{n_{2}}a_{i_{1}i_{2}}%
\overline{a_{j_{1}i_{2}}}\right)  ^{\frac{1}{2}},
\]
by (\ref{delt}), we get
\begin{align*}
\left\vert A\left(  x^{(1)},x^{(2)}\right)  \right\vert  &  \leq\left(
\sum_{\substack{i_{1}=1\\j_{1}=1}}^{n_{2}}y_{i_{1}}^{1}\overline{y_{j_{1}%
}^{(1)}}n_{2}\delta_{i_{1}j_{1}}\right)  ^{\frac{1}{2}}\\
&  =n_{2}^{1/2}\left(  \sum_{i_{1}=1}^{n_{1}}\left\vert x_{i_{1}}%
^{1}\right\vert ^{2}\right)  ^{\frac{1}{2}}\\
&  \leq n_{2}^{1/2}n_{1}^{1/2}%
\end{align*}
and hence%
\[
\left\Vert A\right\Vert \leq n_{2}^{1/2}n_{1}^{1/2}.
\]
Thus
\[
\frac{G_{n_{1}n_{2}}^{(2)}}{\sqrt{n_{1}}\sqrt{n_{2}}}\leq1.
\]
To prove the lower estimate, let us prove that for all bilinear forms
$A:\ell_{\infty}^{n_{1}}\times\ell_{2}^{n_{2}}\longrightarrow\mathbb{C}$ with
unimodular coefficients we have%
\[
\left\Vert A\right\Vert \geq n_{2}^{1/2}n_{1}^{1/2}.
\]
In fact, by the orthogonality of Rademacher functions $r_{i}%
:[0,1]\longrightarrow\{-1,1\}$, we have%
\begin{align*}
n_{2}^{1/2}n_{1}^{1/2}  &  =\left(
{\textstyle\sum\limits_{j=1}^{n_{2}}}
{\textstyle\sum\limits_{i=1}^{n_{1}}}
\left\vert A\left(  e_{i},e_{j}\right)  \right\vert ^{2}\right)  ^{1/2}\\
&  =\left(
{\textstyle\sum\limits_{j=1}^{n_{2}}}
{\textstyle\int\limits_{0}^{1}}
\left\vert
{\textstyle\sum\limits_{i=1}^{n_{1}}}
r_{i}(t)A\left(  e_{i},e_{j}\right)  \right\vert ^{2}dt\right)  ^{1/2}\\
&  =\left(
{\textstyle\int\limits_{0}^{1}}
{\textstyle\sum\limits_{j=1}^{n_{2}}}
\left\vert A\left(
{\textstyle\sum\limits_{i=1}^{n_{1}}}
r_{i}(t)e_{i},e_{j}\right)  \right\vert ^{2}dt\right)  ^{1/2}\\
&  \leq\sup_{t\in\lbrack0,1]}\left(
{\textstyle\sum\limits_{j=1}^{n_{2}}}
\left\vert A\left(
{\textstyle\sum\limits_{i=1}^{n_{1}}}
r_{i}(t)e_{i},e_{j}\right)  \right\vert ^{2}\right)  ^{1/2}\\
&  \leq\left\Vert A\right\Vert \sup_{\left\Vert \varphi\right\Vert _{\ell_{2}%
}=1}\left(
{\textstyle\sum\limits_{j=1}^{n_{2}}}
\left\vert \varphi(e_{j})\right\vert ^{2}\right)  ^{1/2}\\
&  =\left\Vert A\right\Vert .
\end{align*}

\begin{remark} 
		When $n_{1}>n_{2}$ we were not able to design a constructive,
		definitive answer to the problem. Even the correct asymptotic growth behavior of $G^{(2)}_{n_1n_2}$ in this case seems an interesting theoretical question. 
\end{remark}

\section{The game in higher dimensions}

\label{Section Anisotropic}

The Gale--Berlekamp switching game has a natural extension to higher
dimensions. 
Let $m\geq2$ be an integer and let an $n\times\cdots\times n$ array $\left(
a_{j_{1}\cdots j_{m}}\right)  $ of lights be given each either on
($a_{j_{1}\cdots j_{m}}=1$) or off ($a_{j_{1}\cdots j_{m}}=-1$). Let us also
suppose that for each $j_{k}$ there is a switch $x_{j_{k}}^{(k)}$ so that if
the switch is pulled ($x_{j_{k}}^{(k)}=-1$) all of the corresponding lights
are   switched: on to off or off to on. The
goal is to maximize the difference between the lights on and off. As in the
two-dimensional case, maximize the difference between on-lights on and
off-lights is equivalent to estimate%
\[
\max_{x_{j_{1}}^{(1)},\ldots,x_{j_{m}}^{(m)}\in\{-1,1\}}\left\vert
\sum\limits_{j_{1},\ldots,j_{m}=1}^{n}a_{j_{1}\ldots j_{m}}x_{j_{1}}%
^{(1)}\cdots x_{j_{m}}^{(m)}\right\vert
\]
and the extremal problem consists of estimating%
\[
S_{n}=\min\left\{  \max_{x_{j_{1}}^{(1)},\ldots,x_{j_{m}}^{(m)}\in
\{-1,1\}}\left\vert \sum\limits_{j_{1},\ldots,j_{m}=1}^{n}a_{j_{1}\ldots
j_{m}}x_{j_{1}}^{(1)}\cdots x_{j_{m}}^{(m)}\right\vert :a_{j_{1}\ldots j_{m}%
}=1\text{ or }-1\right\}  ,
\]
As in the bilinear case,%
\[
S_{n}=\min\left\Vert A:\ell_{\infty}^{n}\times\cdots\times\ell_{\infty}%
^{n}\rightarrow\mathbb{R}\right\Vert ,
\]
with%
\[
A\left(  x^{(1)},\ldots,x^{(m)}\right)  =\sum\limits_{j_{1},\ldots,j_{m}
=1}^{n}a_{j_{1}\ldots j_{m}}x_{j_{1}}^{(1)}\cdots x_{j_{m}}^{(m)}.
\]
The anisotropic case allows for $n_{1}\times\cdots\times n_{m}$ arrays not necessarily
square arrays and, in this case, we write%
\[
S_{n_{1}\ldots n_{m}}=\min\left\{  \max_{x_{j_{1}}^{(1)},\ldots,x_{j_{m}
}^{(m)}\in\{-1,1\}}\left\vert \sum\limits_{j_{1},\ldots,j_{m}=1}^{n_{1}
,\ldots,n_{m}}a_{j_{1}\ldots j_{m}}x_{j_{1}}^{(1)}\cdots x_{j_{m}}%
^{(m)}\right\vert :a_{j_{1}\ldots j_{m}}=1\text{ or }-1\right\}  .
\]
From a recent result of \cite{al}, we can easily obtain 
\begin{equation}
\frac{1}{m\left(  \sqrt{2}\right)  ^{m-1}}\leq\frac{S_{n_{1}\ldots n_{m}}%
}{\sqrt{n_{1\cdots}n_{m}}\max\{\sqrt{n_{1}},\ldots,\sqrt{n_{m}}\}}\leq
8m\sqrt{m!}\sqrt{\log(1+4m)}. \label{25}%
\end{equation}

Following the notation of \cite{ap5}, let $m\geq2$ be an integer and $\left(
a_{i_{1}\cdots i_{m}}\right)  $ and $n\times\cdots\times n$ be an array of
complex scalars such that $\left\vert a_{i_{1}\cdots i_{m}}\right\vert =1$.
For $p\in(1,\infty]$, let
\[
g_{m,n}^{\mathbb{C}}(p)=\max\left\vert \sum\limits_{i_{1},\ldots,i_{m}=1}%
^{n}a_{i_{1}\cdots i_{m}}x_{i_{1}}^{(1)}\cdots x_{i_{m}}^{(m)}\right\vert ,
\]
where the maximum is evaluated over all $x_{i_{j}}^{(j)}\in\mathbb{C}$ such
that $\Vert(x_{i_{j}}^{(j)})_{i=1}^{n}\Vert_{p}=1$ for all $j.$ It is not
difficult to prove that
\[
g_{m,n}^{\mathbb{C}}(p)=\left\Vert A:\ell_{p}^{n}\times\cdots\times\ell
_{p}^{n}\rightarrow\mathbb{C}\right\Vert ,
\]
with%
\[
A\left( x^{(1)},\ldots,x^{(m)}\right) =\sum\limits_{i_{1},\ldots,i_{m}=1}%
^{n}a_{i_{1}\cdots i_{m}}x_{i_{1}}^{(1)}\cdots x_{i_{m}}^{(m)}.
\]
Denoting
\[
G_{m,n}(p)=\min g_{m,n}^{\mathbb{C}}(p),
\]
where minimum is evaluated over all unimodular $m$-linear forms $A:\ell
_{p}^{n}\times\cdots\times\ell_{p}^{n}\rightarrow\mathbb{C}$, the best
information we can collect (combining results from \cite{ap5, djd}) is the following:%

\begin{equation}
\left\vert
\begin{array}
[c]{l}%
\frac{1}{1.3m^{0.365}}\leq\frac{G_{m,n}(p)}{n^{\frac{mp+p-2m}{2p}}}\leq
8\sqrt{m!}\sqrt{\log(1+4m)}\text{ for }p\in\lbrack2,\infty]\\
1\leq\frac{G_{m,n}(p)}{n^{1-\frac{1}{p}}}\leq C_{m}\text{ for }[1,2],
\end{array}
\right.  \label{7ttl}%
\end{equation}
where $C_{m}$ is obtained by interpolation (via the Riesz--Thorin Theorem) of
the constant $1$ (the constant when $p=1$) and $8\sqrt{m!}\sqrt{\log(1+4m)}$
(the constant when $p=2$). Moreover, defining
\[
g_{m,n}^{\mathbb{C}}(p_{1},\ldots,p_{m})=\left\Vert A:\ell_{p_{1}}^{n}%
\times\cdots\times\ell_{p_{m}}^{n}\rightarrow\mathbb{C}\right\Vert ,
\]
and%
\[
G_{m,n}(p_{1},\ldots,p_{m})=\min g_{m,n}^{\mathbb{C}}(p_{1},\ldots,p_{m}),
\]
from \cite{djd}, it is possible to show that%

\[
K_{m}\leq\frac{G_{m,n}(p_{1},\ldots,p_{m})}{n^{\frac{1}{\min\left\{
\max\{2,p_{k}^{\ast}\right\}  \}}+\sum_{k=1}^{m}\max\left\{  \frac{1}{2}
-\frac{1}{p_{k}},0\right\}  }}\leq D_{m},
\]
where $D_{m},K_{m}$ behave essentially as the constants from (\ref{7ttl}).

The above solution rests in a non-deterministic tool. We
shall show in what follows that for the case of complex scalars we can find deterministic
solutions with better constants, which happen to be optimal in same meaningful cases.

\subsection{{The anisotropic, continuous game}}

We begin with a matrix $\left(  a_{j_{1}\ldots j_{m}}\right)  _{j_{1}
,\ldots,j_{m}=1}^{n_{1},\ldots,n_{m}}$ whose elements are unit vectors in the
Euclidean space $\mathbb{R}^{2}.$ The initial direction pattern of each
$n_{1}\cdots n_{m}$ vectors is set up at the beginning of the game. For each
$k\in\{1,\ldots,m\}$, we have $n_{k}$ knobs $x_{1}^{(k)},\ldots,x_{n_{k}%
}^{(k)}$. When the knob $x_{j_{k}}^{(k)}$ is rotated by an angle
$\theta_{j_{k}}^{(k)},$ the same happens with all the vectors $a_{j_{1}\ldots
j_{m}}$ with $j_{k}$ fixed. Defining $\Theta$ and $s(\Theta)$ as in the
two-dimensional case, the extremal problem is to determine
\[
G_{n_{1}\ldots n_{m}}^{(1)}:=\min\{s(\Theta):\Theta\text{ an }n_{1}%
\times\cdots\times n_{m}\text{ pattern}\}.
\]
We prove the following:

\begin{theorem}
\label{8b}For all positive integers $m\geq2$, $n_{1},\ldots,n_{m}\geq1$ we
have
\[
(0.886)^{m-1}\simeq\left(  \frac{\sqrt{\pi}}{2}\right)  ^{m-1}\leq
\frac{G_{n_{1}\ldots n_{m}}^{(1)}}{\sqrt{n_{1}\cdots n_{m}}\max\{\sqrt{n_{1}
},\ldots,\sqrt{n_{m}}\}}\leq1.
\]
Moreover, the universal upper bound $1$ cannot be improved.
\end{theorem}

The proof that, in general, the upper bound $1$ cannot be improved is trivial ---
 just consider $n_{2}=\cdots=n_{m}=1$ and note that in this case%
\[
G_{n_{1}\ldots n_{m}}^{(1)}=n_{1}=\sqrt{n_{1}\cdots n_{m}}\max\{\sqrt{n_{1}%
},\ldots,\sqrt{n_{m}}\}.
\]
The case $n_{1}=\cdots=n_{m}$ was investigated in \cite{ap5}, but the
techniques from \cite{ap5} do not provide good estimates for the upper
constants: for instance, if we follow the arguments from \cite{ap5} we just
obtain $8\sqrt{m!}\sqrt{\log(1+4m)}$, due to Kahane--Salem--Zygmund
inequality, instead of the universal sharp constant $1$.

\bigskip

\subsubsection{Proof of Theorem \ref{8b}}

Our task is to estimate $\inf\left\{  \left\Vert A\right\Vert :\left\vert
a_{j_{1}\ldots j_{m}}\right\vert =1\right\}  $, where the infimum runs over
all $m$-linear forms $A:\ell_{\infty}^{n_{1}}\times\cdots\times\ell_{\infty
}^{n_{m}}\rightarrow\mathbb{C}$ with unimodular coefficients.

With no loss of generality, we suppose  $n_{1}\leq\cdots\leq n_{m}$.
For upper bound, consider, for all $k=1,\ldots,m-1$, a $n_{k+1}\times n_{k+1}$
matrix $\left( a_{rs}^{(k)}\right) $ with%
\[
a_{rs}^{(k)}=e^{2\pi i\frac{rs}{n_{k+1}}}.
\]
A simple computation shows that
\begin{align*}
&
\begin{cases}
\sum_{t=1}^{n_{2}}a_{rt}^{(1)}\overline{a_{st}^{(1)}}=n_{2}\delta_{rs}.\\
|a_{rs}^{(1)}|=1.
\end{cases}
\\
&  \vdots\\
&
\begin{cases}
\sum_{t=1}^{n_{m}}a_{rt}^{(m-1)}\overline{a_{st}^{(m-1)}}=n_{m}\delta_{rs}.\\
|a_{rs}^{(m-1)}|=1.
\end{cases}
\end{align*}
All the matrices are completed with zeros (if necessary) in order to get a
square matrix $n_{m}\times n_{m}.$ Define%
\[
A:\ell_{\infty}^{n_{1}}\times\cdots\times\ell_{\infty}^{n_{m}}\rightarrow
\mathbb{C}%
\]
by
\[
A\left(  x^{(1)},\dots,x^{(m)}\right)  =\sum_{i_{1},\dots,i_{m}=1}%
^{n_{1},\ldots,n_{m}}a_{i_{1}i_{2}}^{(1)}a_{i_{2}i_{3}}^{(2)}\cdots
a_{i_{m-1}i_{m}}^{(m-1)}x_{i_{1}}^{(1)}\cdots x_{i_{m}}^{(m)}%
\]
and note that, since $n_{1}\leq\cdots\leq n_{m}$, the coefficients
\[
c_{i_{1}\cdots i_{m}}:=a_{i_{1}i_{2}}^{(1)}a_{i_{2}i_{3}}^{(2)}\cdots
a_{i_{m-1}i_{m}}^{(m-1)}%
\]
of all monomials $x_{i_{1}}^{(1)}\cdots x_{i_{m}}^{(m)}$ with $i_{k}%
\in\{1,\ldots,n_{k}\}$ are unimodular. For $x^{(1)}\in B_{\ell_{\infty}
^{n_{1}}},\ldots,x^{(m)}\in B_{\ell_{\infty}^{n_{m}}},$ consider $y^{(1)}\in
B_{\ell_{\infty}^{n_{m}}},\ldots,y^{(m)}\in B_{\ell_{{\infty}}^{n_{m}}}$
defined by
\[
y^{(1)}=\left( x_{1}^{(1)},\ldots,x_{n_{1}}^{(1)},0,\ldots,0\right)
\]
and so on. We have%
\begin{align*}
&  \left\vert A\left(  x^{(1)},\dots,x^{(m)}\right)  \right\vert \\
&  =\left\vert \sum_{i_{1},\dots,i_{m}=1}^{n_{m}}a_{i_{1}i_{2}}^{(1)}%
a_{i_{2}i_{3}}^{(2)}\cdots a_{i_{m-1}i_{m}}^{(m-1)}y_{i_{1}}^{(1)}\cdots
y_{i_{m}}^{(m)}\right\vert \\
&  \leq\sum_{i_{m}=1}^{n_{m}}\left\vert \sum_{i_{1},\dots,i_{m-1}=1}^{n_{m}%
}a_{i_{1}i_{2}}^{(1)}a_{i_{2}i_{3}}^{(2)}\cdots a_{i_{m-1}i_{m}}%
^{(m-1)}y_{i_{1}}^{(1)}\cdots y_{i_{m-1}}^{(m-1)}\right\vert |y_{i_{m}}^{m}|\\
&  \leq\left(  \sum_{i_{m}=1}^{n_{m}}|y_{i_{m}}^{(m)}|^{2}\right)  ^{1/2}%
\cdot\left(  \sum_{i_{m}=1}^{n_{m}}\left\vert \sum_{i_{1},\dots,i_{m-1}
=1}^{n_{m}}a_{i_{1}i_{2}}^{(1)}a_{i_{2}i_{3}}^{(2)}\cdots a_{i_{m-1}i_{m}%
}^{(m-1)}y_{i_{1}}^{(1)}\cdots y_{i_{m-1}}^{(m-1)}\right\vert ^{2}\right)
^{1/2}\\
&  \leq n_{m}^{1/2}\left(  \sum_{i_{m}=1}^{n_{m}}\sum_{\substack{i_{1}
,\dots,i_{m-1}=1\\j_{1},\dots,j_{m-1}=1}}^{n_{m}}a_{i_{1}i_{2}}^{(1)}%
\overline{a_{j_{1}j_{2}}^{(1)}}\cdots a_{i_{m-1}i_{m}}^{(m-1)}\overline
{a_{j_{m-1}i_{m}}^{(m-1)}}y_{i_{1}}^{(1)}\overline{y_{j_{1}}^{(1)}}\cdots
y_{i_{m-1}}^{(m-1)}\overline{y_{j_{m-1}}^{(m-1)}}\right)  ^{1/2}.
\end{align*}
Thus%
\begin{align*}
&  \left\vert A\left(  x^{(1)},\dots,x^{(m)}\right)  \right\vert \\
&  \leq n_{m}^{1/2}\left(  \sum_{i_{m}=1}^{n_{m}}\sum_{\substack{i_{1}
,\dots,i_{m-1}=1\\j_{1},\dots,j_{m-1}=1}}^{n_{m}}a_{i_{1}i_{2}}^{(1)}%
\overline{a_{j_{1}j_{2}}^{(1)}}\cdots a_{i_{m-1}i_{m}}^{(m-1)}\overline
{a_{j_{m-1}i_{m}}^{(m-1)}}y_{i_{1}}^{(1)}\overline{y_{j_{1}}^{(1)}}\cdots
y_{i_{m-1}}^{(m-1)}\overline{y_{j_{m-1}}^{(m-1)}}\right)  ^{1/2}\\
&  \leq n_{m}^{1/2}\left(  \sum_{\substack{i_{1},\dots,i_{m-1}=1\\j_{1}
,\dots,j_{m-1}=1}}^{n_{m}}a_{i_{1}i_{2}}^{(1)}\overline{a_{j_{1}j_{2}}^{(1)}%
}\cdots a_{i_{m-2}i_{m-1}}^{(m-2)}\overline{a_{j_{m-2}j_{m-1}}^{(m-2)}%
}y_{i_{1}}^{(1)}\overline{y_{j_{1}}^{(1)}}\cdots y_{i_{m-1}}^{(m-1)}%
\overline{y_{j_{m-1}}^{(m-1)}}\sum_{i_{m}=1}^{n_{m}}a_{i_{m-1}i_{m}}%
^{(m-1)}\overline{a_{j_{m-1}i_{m}}^{(m-1)}}\right)  ^{1/2}.
\end{align*}
Since
\[
\sum_{i_{m}=1}^{n_{m}}a_{i_{m-1}i_{m}}^{(m-1)}\overline{a_{j_{m-1}i_{m}
}^{(m-1)}}=n_{m}\delta_{i_{m-1}j_{m-1}},
\]
we have%
\begin{align*}
&  \left\vert A\left(  x^{(1)},\dots,x^{(m)}\right)  \right\vert \\
&  \leq n_{m}^{1/2}\left(  \sum_{i_{m-1}=1}^{n_{m}}\sum_{\substack{i_{1}
,\dots,i_{m-2}=1\\j_{1},\dots,j_{m-2}=1}}^{n_{m}}a_{i_{1}i_{2}}^{(1)}%
\overline{a_{j_{1}j_{2}}^{(1)}}\cdots a_{i_{m-2}i_{m-1}}^{(m-2)}%
\overline{a_{j_{m-2}i_{m-1}}^{(m-2)}}y_{i_{1}}^{(1)}\overline{y_{j_{1}}^{(1)}%
}\cdots y_{i_{m-2}}^{(m-2)}\overline{y_{j_{m-2}}^{(m-2)}}y_{i_{m-1}}%
^{(m-1)}\overline{y_{i_{m-1}}^{(m-1)}}n_{m}\right)  ^{1/2}\\
&  =n_{m}\left(  \sum_{i_{m-1}=1}^{n_{m}}\sum_{\substack{i_{1},\dots
,i_{m-2}=1\\j_{1},\dots,j_{m-2}=1}}^{n_{m}}a_{i_{1}i_{2}}^{(1)}\overline
{a_{j_{1}j_{2}}^{(1)}}\cdots a_{i_{m-2}i_{m-1}}^{(m-2)}\overline
{a_{j_{m-2}i_{m-1}}^{(m-2)}}y_{i_{1}}^{(1)}\overline{y_{j_{1}}^{(1)}}\cdots
y_{i_{m-2}}^{(m-2)}\overline{y_{j_{m-2}}^{(m-2)}}\left\vert y_{i_{m-1}%
}^{(m-1)}\right\vert ^{2}\right)  ^{1/2}\\
&  \leq n_{m}\left(  \sum_{\substack{i_{1},\dots,i_{m-2}=1\\j_{1}
,\dots,j_{m-2}=1}}^{n_{m}}a_{i_{1}i_{2}}^{(1)}\overline{a_{j_{1}j_{2}}^{(1)}%
}\cdots a_{i_{m-3}i_{m-2}}^{(m-3)}\overline{a_{j_{m-3}i_{m-2}}^{(m-3)}%
}y_{i_{1}}^{(1)}\overline{y_{j_{1}}^{(1)}}\cdots y_{i_{m-2}}^{(m-2)}%
\overline{y_{j_{m-2}}^{(m-2)}}\sum_{i_{m-1}=1}^{n_{m}}a_{i_{m-2}i_{m-1}%
}^{(m-2)}\overline{a_{j_{m-2}i_{m-1}}^{(m-2)}}\right)  ^{1/2}.
\end{align*}
Thus
\begin{align*}
&  \left\vert A\left(  x^{(1)},\dots,x^{(m)}\right)  \right\vert \\
&  \leq n_{m}\left(  \sum_{\substack{i_{1},\dots,i_{m-2}=1\\j_{1}
,\dots,j_{m-2}=1}}^{n_{m}}a_{i_{1}i_{2}}^{(1)}\overline{a_{j_{1}j_{2}}^{(1)}%
}\cdots a_{i_{m-3}i_{m-2}}^{(m-3)}\overline{a_{j_{m-3}i_{m-2}}^{(m-3)}%
}y_{i_{1}}^{(1)}\overline{y_{j_{1}}^{(1)}}\cdots y_{i_{m-2}}^{(m-2)}%
\overline{y_{j_{m-2}}^{(m-2)}}\sum_{i_{m-1}=1}^{n_{m}}a_{i_{m-2}i_{m-1}%
}^{(m-2)}\overline{a_{j_{m-2}i_{m-1}}^{(m-2)}}\right)  ^{1/2}.
\end{align*}
Since%
\begin{align*}
&  \left(  \sum_{\substack{i_{1},\dots,i_{m-2}=1\\j_{1},\dots,j_{m-2}
=1}}^{n_{m}}a_{i_{1}i_{2}}^{(1)}\overline{a_{j_{1}j_{2}}^{(1)}}\cdots
a_{i_{m-3}i_{m-2}}^{(m-3)}\overline{a_{j_{m-3}i_{m-2}}^{(m-3)}}y_{i_{1}}%
^{(1)}\overline{y_{j_{1}}^{(1)}}\cdots y_{i_{m-2}}^{(m-2)}\overline
{y_{j_{m-2}}^{(m-2)}}\sum_{i_{m-1}=1}^{n_{m}}a_{i_{m-2}i_{m-1}}^{(m-2)}%
\overline{a_{j_{m-2}i_{m-1}}^{(m-2)}}\right)  ^{1/2}\\
&  =n_{m-1}^{1/2}\left(  \sum_{i_{m-2}=1}^{n_{m}}\sum_{\substack{i_{1}
,\dots,i_{m-3}=1\\j_{1},\dots,j_{m-3}=1}}^{n_{m}}a_{i_{1}i_{2}}^{(1)}%
\overline{a_{j_{1}j_{2}}^{(1)}}\cdots a_{i_{m-3}i_{m-2}}^{(m-3)}%
\overline{a_{j_{m-3}i_{m-2}}^{(m-3)}}y_{i_{1}}^{(1)}\overline{y_{j_{1}}^{(1)}%
}\cdots y_{i_{m-2}}^{(m-2)}\overline{y_{i_{m-2}}^{(m-2)}}\right)  ^{1/2}\\
&  =n_{m-1}^{1/2}\left(  \sum_{i_{m-2}=1}^{n_{m}}|y_{i_{m-2}}^{(m-2)}|^{2}%
\sum_{\substack{i_{1},\dots,i_{m-3}=1\\j_{1},\dots,j_{m-3}=1}}^{n_{m}}%
a_{i_{1}i_{2}}^{(1)}\overline{a_{j_{1}j_{2}}^{(1)}}\cdots a_{i_{m-3}i_{m-2}%
}^{(m-3)}\overline{a_{j_{m-3}i_{m-2}}^{(m-3)}}y_{i_{1}}^{(1)}\overline
{y_{j_{1}}^{(1)}}\cdots y_{i_{m-3}}^{(m-3)}\overline{y_{j_{m-3}}^{(m-3)}%
}\right)  ^{1/2},
\end{align*}
we conclude that%
\begin{align*}
&  \left\vert A\left(  x^{(1)},\dots,x^{(m)}\right)  \right\vert \\
&  \leq n_{m}n_{m-1}^{1/2}\left(  \sum_{i_{m-2}=1}^{n_{m}}\sum
_{\substack{i_{1},\dots,i_{m-3}=1\\j_{1},\dots,j_{m-3}=1}}^{n_{m}}%
a_{i_{1}i_{2}}^{(1)}\overline{a_{j_{1}j_{2}}^{(1)}}\cdots a_{i_{m-3}i_{m-2}%
}^{(m-3)}\overline{a_{j_{m-3}i_{m-2}}^{(m-3)}}y_{i_{1}}^{(1)}\overline
{y_{j_{1}}^{(1)}}\cdots y_{i_{m-3}}^{(m-3)}\overline{y_{j_{m-3}}^{(m-3)}%
}\right)  ^{1/2}%
\end{align*}
and repeating this procedure we finally obtain%
\begin{align*}
\left\vert A\left(  x^{(1)},\dots,x^{(m)}\right)  \right\vert  &  \leq
n_{m}n_{m-1}^{\frac{1}{2}}\cdots n_{2}^{\frac{1}{2}}\left(  \sum_{i_{1}
=1}^{n_{m}}y_{i_{1}}^{(1)}\overline{y_{i_{1}}^{(1)}}\right)  ^{\frac{1}{2}}\\
&  =n_{m}^{\frac{1}{2}}n_{m}^{\frac{1}{2}}\cdots n_{2}^{\frac{1}{2}}\left(
\sum_{i_{1}=1}^{n_{1}}|x_{i_{1}}^{(1)}|^{2}\right)  ^{\frac{1}{2}}\\
&  \leq n_{m}^{\frac{1}{2}}\left(  n_{m}^{\frac{1}{2}}\cdots n_{1}^{\frac
{1}{2}}\right)  .
\end{align*}
Thus
\[
\frac{G_{n_{1}\ldots n_{m}}^{(1)}}{\sqrt{n_{1}\cdots n_{m}}\max\{\sqrt{n_{1}
},\ldots,\sqrt{n_{m}}\}}\leq1.
\]

 The lower estimate is a adaptation of the bilinear case,  {using this well--know extension of inequality (5) to multiple sums as follow:}
 
 \begin{align*}
 &  \left(  \sum\limits_{j_{1},\ldots,j_{m}=1}^{n_{1},\ldots,n_{m}}\left\vert
 a_{j_{1}\ldots j_{m}}\right\vert ^{2}\right)  ^{1/2}\label{est3}\\
 &  \leq\left(  \frac{2}{\sqrt{\pi}}\right)  ^{m}\left(  \frac{1}{2\pi}\right)
 ^{n_{1}+\cdots+n_{m}}\int_{0}^{2\pi}\ldots\int_{0}^{2\pi}\left\vert
 \sum_{j_{1},\ldots,j_{m}=1}^{n_{1},\ldots,n_{m}}a_{j_{1}\ldots j_{m}%
 }e^{it_{j_{1}}^{\left(  1\right)  }}\cdots e^{it_{j_{m}}^{\left(  m\right)  }%
 }\right\vert dt_{j_{1}}^{\left(  1\right)  }\cdots dt_{j_{m}}^{\left(
 	m\right)  }.
 \end{align*}

\subsection{The weighted, anisotropic continuous game}

We begin with a matrix $\left(  a_{j_{1}\ldots j_{m}}\right)  _{j_{1}
,\ldots,j_{m}=1}^{n_{1},\ldots,n_{m}}$ whose elements are unit vectors in the
Euclidean space $\mathbb{R}^{2}.$ The initial direction pattern of each
$n_{1}\cdots n_{m}$ vectors is set up at the beginning of the game. For each
$k\in\{1,\ldots,m\}$, we have $n_{k}$ knobs $x_{1}^{(k)},\ldots,x_{n_{k}%
}^{(k)}$. When the knob $x_{j_{k}}^{(k)}$ is rotated by an angle
$\theta_{j_{k}}^{(k)},$ the same happens with all the vectors $a_{j_{1}\ldots
j_{m}}$ with $j_{k}$ fixed. However, in the case $k=m$, the respective vectors
$a_{{j_{1}\ldots j_{m}}}$, in addition to being rotated, are also multiplied
by a real number $b_{j_{m}},$ with%
\[%
{\textstyle\sum\limits_{j_{m}=1}^{n_{m}}}
b_{j_{m}}^{2}=1.
\]
Defining $\Theta$ and $s(\Theta)$ as in the two-dimensional case, the extremal
problem is to determine
\[
G_{n_{1}\ldots n_{m}}^{(2)}:=\min\{s(\Theta):\Theta\text{ an }n_{1}%
\times\cdots\times n_{m}\text{ pattern}\}.
\]
We prove the following:

\begin{theorem}
\label{9b}For all positive integers $m,n_{1},\ldots,n_{m}\geq2$, with
$n_{m}\geq\max\{n_{1},\ldots,n_{m-1}\}$, we have
\[
\frac{G_{n_{1}\ldots n_{m}}^{(2)}}{\sqrt{n_{1}}\cdots\sqrt{n_{m}}}=1.
\]

\end{theorem}

\subsubsection{Proof of Theorem \ref{9b}}

With no loss of generality,  suppose that $n_{1}\leq\cdots\leq n_{m}$.
For the upper estimate, consider, in as Theorem \ref{8b}, for all
$k=1,\ldots,m-1$, a $n_{k+1}\times n_{k+1}$ matrix $\left( a_{rs}^{(k)}\right)
$ with%
\[
a_{rs}^{(k)}=e^{2\pi i\frac{rs}{n_{k+1}}}.
\]
Then, the $m$-linear form%
\[
A:\ell_{\infty}^{n_{1}}\times\cdots\times\ell_{\infty}^{n_{m-1}}\times\ell
_{2}^{n_{m}}\rightarrow\mathbb{C}%
\]
defined by
\[
A\left(  x^{(1)},\dots,x^{(m)}\right)  =\sum_{i_{1},\dots,i_{m}=1}%
^{n_{1},\ldots,n_{m}}a_{i_{1}i_{2}}^{(1)}a_{i_{2}i_{3}}^{(2)}\cdots
a_{i_{m-1}i_{m}}^{(m-1)}x_{i_{1}}^{(1)}\cdots x_{i_{m}}^{(m)},
\]
is unimodular. For $x^{(1)}\in B_{\ell_{\infty}^{n_{1}}},\ldots,x^{(m-1)}\in
B_{\ell_{\infty}^{n_{m-1}}},x^{(m)}\in B_{\ell_{2}^{n_{m}}}$ consider
$y^{(1)}\in B_{\ell_{{\infty}}^{n_{m}}},\ldots,y^{(m-1)}\in B_{\ell_{{\infty}
}^{n_{m}}},y^{(m)}\in B_{\ell_{2}^{n_{m}}}$ defined by
\[
y^{(1)}=\left( x_{1}^{(1)},\ldots,x_{n_{1}}^{(1)},0,\ldots,0\right)
\]
and so on. We have%
\begin{align*}
&  \left\vert A\left(  x^{(1)},\dots,x^{(m)}\right)  \right\vert \\
&  =\left\vert \sum_{i_{1},\dots,i_{m}=1}^{n_{m}}a_{i_{1}i_{2}}^{(1)}%
a_{i_{2}i_{3}}^{(2)}\cdots a_{i_{m-1}i_{m}}^{(m-1)}y_{i_{1}}^{(1)}\cdots
y_{i_{m}}^{(m)}\right\vert \\
&  \leq\sum_{i_{m}=1}^{n_{m}}\left\vert \sum_{i_{1},\dots,i_{m-1}=1}^{n_{m}%
}a_{i_{1}i_{2}}^{(1)}a_{i_{2}i_{3}}^{(2)}\cdots a_{i_{m-1}i_{m}}%
^{(m-1)}y_{i_{1}}^{(1)}\cdots y_{i_{m-1}}^{(m-1)}\right\vert |y_{i_{m}}^{m}|\\
&  \leq\left(  \sum_{i_{m}=1}^{n_{m}}|y_{i_{m}}^{(m)}|^{2}\right)  ^{1/2}%
\cdot\left(  \sum_{i_{m}=1}^{n_{m}}\left\vert \sum_{i_{1},\dots,i_{m-1}
=1}^{n_{m}}a_{i_{1}i_{2}}^{(1)}a_{i_{2}i_{3}}^{(2)}\cdots a_{i_{m-1}i_{m}%
}^{(m-1)}y_{i_{1}}^{(1)}\cdots y_{i_{m-1}}^{(m-1)}\right\vert ^{2}\right)
^{1/2}\\
&  \leq\left(  \sum_{i_{m}=1}^{n_{m}}\left\vert \sum_{i_{1},\dots,i_{m-1}
=1}^{n_{m}}a_{i_{1}i_{2}}^{(1)}a_{i_{2}i_{3}}^{(2)}\cdots a_{i_{m-1}i_{m}%
}^{(m-1)}y_{i_{1}}^{(1)}\cdots y_{i_{m-1}}^{(m-1)}\right\vert ^{2}\right)
^{1/2}.
\end{align*}
Since%
\begin{align*}
&  \left(  \sum_{i_{m}=1}^{n_{m}}\left\vert \sum_{i_{1},\dots,i_{m-1}
=1}^{n_{m}}a_{i_{1}i_{2}}^{(1)}a_{i_{2}i_{3}}^{(2)}\cdots a_{i_{m-1}i_{m}%
}^{(m-1)}y_{i_{1}}^{(1)}\cdots y_{i_{m-1}}^{(m-1)}\right\vert ^{2}\right)
^{1/2}\\
&  =\left(  \sum_{i_{m}=1}^{n_{m}}\sum_{\substack{i_{1} ,\dots,i_{m-1}
=1\\j_{1},\dots,j_{m-1}=1}}^{n_{m}}a_{i_{1}i_{2}}^{(1)}\overline{a_{j_{1}
j_{2}}^{(1)}}\cdots a_{i_{m-1}i_{m}}^{(m-1)}\overline{a_{j_{m-1}i_{m}}
^{(m-1)}}y_{i_{1}}^{(1)}\overline{y_{j_{1}}^{(1)}}\cdots y_{i_{m-1}}%
^{(m-1)}\overline{y_{j_{m-1}}^{(m-1)}}\right)  ^{1/2}\\
&  =\left(  \sum_{\substack{i_{1},\dots,i_{m-1}=1\\j_{1} ,\dots,j_{m-1}
=1}}^{n_{m}}a_{i_{1}i_{2}}^{(1)}\overline{a_{j_{1}j_{2}}^{(1)}}\cdots
a_{i_{m-2}i_{m-1}}^{(m-2)}\overline{a_{j_{m-2}j_{m-1}}^{(m-2)}}y_{i_{1}}%
^{(1)}\overline{y_{j_{1}}^{(1)}}\cdots y_{i_{m-1}}^{(m-1)}\overline
{y_{j_{m-1}}^{(m-1)}}\sum_{i_{m}=1}^{n_{m}}a_{i_{m-1}i_{m}}^{(m-1)}%
\overline{a_{j_{m-1}i_{m}}^{(m-1)}}\right)  ^{1/2}%
\end{align*}
and
\[
\sum_{i_{m}=1}^{n_{m}}a_{i_{m-1}i_{m}}^{(m-1)}\overline{a_{j_{m-1}i_{m}
}^{(m-1)}}=n_{m}\delta_{i_{m-1}j_{m-1}},
\]
we have%
\begin{align*}
&  \left\vert A\left(  x^{(1)},\dots,x^{(m)}\right)  \right\vert \\
&  \leq\left(  \sum_{i_{m-1}=1}^{n_{m}}\sum_{\substack{i_{1} ,\dots
,i_{m-2}=1\\j_{1},\dots,j_{m-2}=1}}^{n_{m}}a_{i_{1}i_{2}}^{(1)}\overline
{a_{j_{1}j_{2}}^{(1)}}\cdots a_{i_{m-2}i_{m-1}}^{(m-2)}\overline
{a_{j_{m-2}i_{m-1}}^{(m-2)}}y_{i_{1}}^{(1)}\overline{y_{j_{1}}^{(1)}}\cdots
y_{i_{m-2}}^{(m-2)}\overline{y_{j_{m-2}}^{(m-2)}}y_{i_{m-1}}^{(m-1)}%
\overline{y_{i_{m-1}}^{(m-1)}}n_{m}\right)  ^{1/2}\\
&  =n_{m}^{1/2}\left(  \sum_{i_{m-1}=1}^{n_{m}}\sum_{\substack{i_{1}
,\dots,i_{m-2}=1\\j_{1},\dots,j_{m-2}=1}}^{n_{m}}a_{i_{1}i_{2}}^{(1)}%
\overline{a_{j_{1}j_{2}}^{(1)}}\cdots a_{i_{m-2}i_{m-1}}^{(m-2)}%
\overline{a_{j_{m-2}i_{m-1}}^{(m-2)}}y_{i_{1}}^{(1)}\overline{y_{j_{1}}^{(1)}%
}\cdots y_{i_{m-2}}^{(m-2)}\overline{y_{j_{m-2}}^{(m-2)}}\left\vert
y_{i_{m-1}}^{(m-1)}\right\vert ^{2}\right)  ^{1/2}\\
&  \leq n_{m}^{1/2}\left(  \sum_{\substack{i_{1},\dots,i_{m-2}=1\\j_{1}
,\dots,j_{m-2}=1}}^{n_{m}}a_{i_{1}i_{2}}^{(1)}\overline{a_{j_{1}j_{2}}^{(1)}%
}\cdots a_{i_{m-3}i_{m-2}}^{(m-3)}\overline{a_{j_{m-3}i_{m-2}}^{(m-3)}%
}y_{i_{1}}^{(1)}\overline{y_{j_{1}}^{(1)}}\cdots y_{i_{m-2}}^{(m-2)}%
\overline{y_{j_{m-2}}^{(m-2)}}\sum_{i_{m-1}=1}^{n_{m}}a_{i_{m-2}i_{m-1}%
}^{(m-2)}\overline{a_{j_{m-2}i_{m-1}}^{(m-2)}}\right)  ^{1/2}.
\end{align*}
Thus
\begin{align*}
&  \left\vert A\left(  x^{(1)},\dots,x^{(m)}\right)  \right\vert \\
&  \leq n_{m}^{1/2}\left(  \sum_{\substack{i_{1},\dots,i_{m-2}=1\\j_{1}
,\dots,j_{m-2}=1}}^{n_{m}}a_{i_{1}i_{2}}^{(1)}\overline{a_{j_{1}j_{2}}^{(1)}%
}\cdots a_{i_{m-3}i_{m-2}}^{(m-3)}\overline{a_{j_{m-3}i_{m-2}}^{(m-3)}%
}y_{i_{1}}^{(1)}\overline{y_{j_{1}}^{(1)}}\cdots y_{i_{m-2}}^{(m-2)}%
\overline{y_{j_{m-2}}^{(m-2)}}\sum_{i_{m-1}=1}^{n_{m}}a_{i_{m-2}i_{m-1}%
}^{(m-2)}\overline{a_{j_{m-2}i_{m-1}}^{(m-2)}}\right)  ^{1/2}.
\end{align*}
Since%
\begin{align*}
&  \left(  \sum_{\substack{i_{1},\dots,i_{m-2}=1\\j_{1},\dots,j_{m-2}
=1}}^{n_{m}}a_{i_{1}i_{2}}^{(1)}\overline{a_{j_{1}j_{2}}^{(1)}}\cdots
a_{i_{m-3}i_{m-2}}^{(m-3)}\overline{a_{j_{m-3}i_{m-2}}^{(m-3)}}y_{i_{1}}%
^{(1)}\overline{y_{j_{1}}^{(1)}}\cdots y_{i_{m-2}}^{(m-2)}\overline
{y_{j_{m-2}}^{(m-2)}}\sum_{i_{m-1}=1}^{n_{m}}a_{i_{m-2}i_{m-1}}^{(m-2)}%
\overline{a_{j_{m-2}i_{m-1}}^{(m-2)}}\right)  ^{1/2}\\
&  =n_{m-1}^{1/2}\left(  \sum_{i_{m-2}=1}^{n_{m}}\sum_{\substack{i_{1}
,\dots,i_{m-3}=1\\j_{1},\dots,j_{m-3}=1}}^{n_{m}}a_{i_{1}i_{2}}^{(1)}%
\overline{a_{j_{1}j_{2}}^{(1)}}\cdots a_{i_{m-3}i_{m-2}}^{(m-3)}%
\overline{a_{j_{m-3}i_{m-2}}^{(m-3)}}y_{i_{1}}^{(1)}\overline{y_{j_{1}}^{(1)}%
}\cdots y_{i_{m-2}}^{(m-2)}\overline{y_{i_{m-2}}^{(m-2)}}\right)  ^{1/2}\\
&  =n_{m-1}^{1/2}\left(  \sum_{i_{m-2}=1}^{n_{m}}|y_{i_{m-2}}^{(m-2)}|^{2}%
\sum_{\substack{i_{1},\dots,i_{m-3}=1\\j_{1},\dots,j_{m-3}=1}}^{n_{m}}%
a_{i_{1}i_{2}}^{(1)}\overline{a_{j_{1}j_{2}}^{(1)}}\cdots a_{i_{m-3}i_{m-2}%
}^{(m-3)}\overline{a_{j_{m-3}i_{m-2}}^{(m-3)}}y_{i_{1}}^{(1)}\overline
{y_{j_{1}}^{(1)}}\cdots y_{i_{m-3}}^{(m-3)}\overline{y_{j_{m-3}}^{(m-3)}%
}\right)  ^{1/2},
\end{align*}
we conclude that%
\begin{align*}
&  \left\vert A\left(  x^{(1)},\dots,x^{(m)}\right)  \right\vert \\
&  \leq n_{m}^{1/2}n_{m-1}^{1/2}\left(  \sum_{i_{m-2}=1}^{n_{m}}%
\sum_{\substack{i_{1},\dots,i_{m-3}=1\\j_{1},\dots,j_{m-3}=1}}^{n_{m}}%
a_{i_{1}i_{2}}^{(1)}\overline{a_{j_{1}j_{2}}^{(1)}}\cdots a_{i_{m-3}i_{m-2}%
}^{(m-3)}\overline{a_{j_{m-3}i_{m-2}}^{(m-3)}}y_{i_{1}}^{(1)}\overline
{y_{j_{1}}^{(1)}}\cdots y_{i_{m-3}}^{(m-3)}\overline{y_{j_{m-3}}^{(m-3)}%
}\right)  ^{1/2}%
\end{align*}
and repeating this procedure we finally obtain%
\begin{align*}
\left\vert A\left(  x^{(1)},\dots,x^{(m)}\right)  \right\vert  &  \leq
n_{m}^{\frac{1}{2}}n_{m-1}^{\frac{1}{2}}\cdots n_{2}^{\frac{1}{2}}\left(
\sum_{i_{1} =1}^{n_{m}}y_{i_{1}}^{(1)}\overline{y_{i_{1}}^{(1)}}\right)
^{\frac{1}{2}}\\
&  =n_{m}^{\frac{1}{2}}\cdots n_{2}^{\frac{1}{2}}\left(  \sum_{i_{1}=1}%
^{n_{1}}|x_{i_{1}}^{(1)}|^{2}\right)  ^{\frac{1}{2}}\\
&  \leq n_{m}^{\frac{1}{2}}\cdots n_{1}^{\frac{1}{2}}.
\end{align*}
Thus%

\[
\frac{G_{n_{1}\ldots n_{m}}^{(2)}}{\sqrt{n_{1}}\cdots\sqrt{n_{m}}}\leq1.
\]

The lower estimate an adaptation of the argument used in the proof of Theorem
\ref{7b}.

\bigskip

\noindent\textbf{Acknowledgments} The authors thank Fernando Costa J\'unior
for kindly designing and providing  Figure 1.

\end{document}